\documentclass{svproc}
\usepackage{times}
\usepackage{soul}
\usepackage{url}

\usepackage{algcompatible,algorithm,amsmath,amssymb}
\usepackage{amsfonts}
\usepackage{graphicx,multicol,footmisc}

\newtheorem{lem}{Lemma}

\newcommand{\mR}{\mathbb{R}}

\newcommand{\mE}{\mathbb{E}}

\newcommand{\mL}{\mathcal{L}}

\newcommand{\vx}{\boldsymbol{x}}
\newcommand{\vv}{\boldsymbol{v}}
\newcommand{\vw}{\boldsymbol{w}}
\newcommand{\vz}{\boldsymbol{z}}
\newcommand{\vu}{\boldsymbol{u}}

\newcommand{\parallelsum}{\mathbin{\!/\mkern-5mu/\!}}
\newcommand\numberthis{\addtocounter{equation}{1}\tag{\theequation}}

\begin{document}
\mainmatter

\title{Convergence of a Relaxed Variable Splitting Method for Learning Sparse Neural Networks via $\ell_1, \ell_0$, and Transformed-$\ell_1$ Penalties}
\titlerunning{Convergence of RVSM}

\author{
 Thu Dinh
 \and
 Jack Xin}
 \institute{Department of Mathematics,\\
University of California, Irvine, CA 92697, USA. \\
 \email{
 \{t.dinh, jack.xin\}@uci.edu
 }}
 \authorrunning{Thu Dinh, Jack Xin}

\maketitle

\begin{abstract}
Sparsification of neural networks is one of the effective complexity reduction methods to improve efficiency and generalizability. We consider the problem of learning a one hidden layer convolutional neural network with ReLU activation function via gradient descent under sparsity promoting penalties. It is known that when the input data is Gaussian distributed, no-overlap networks (without penalties) in regression problems  with ground truth can be learned in polynomial time at high probability. We propose a relaxed variable splitting method integrating thresholding and gradient descent to overcome the non-smoothness in the loss function. The sparsity in network weight is realized during the optimization (training) process. We prove that under $\ell_1, \ell_0,$ and transformed-$\ell_1$ penalties, no-overlap networks can be learned with high probability, and the iterative weights converge to a global limit which is a transformation of the true weight under a novel  thresholding operation. Numerical experiments confirm theoretical findings, and compare the accuracy and sparsity trade-off among the penalties.
\keywords{regularization, sparsification, non-convex optimization}
\end{abstract}

% \medskip

% \hspace {.1 in} {\bf To appear in  Intelligent Systems Conference (IntelliSys) 2020,}

% \hspace{.1 in} {\bf September 3-4,  Amsterdam, The Netherlands.} 

\section{Introduction}
Deep neural networks (DNN) have achieved state-of-the-art performance on many machine learning tasks such as 
speech recognition (Hinton et al., 2012 \cite{Hinton}), computer vision (Krizhevsky et al., 2016 \cite{Krizhevsky}), and 
natural language processing (Dauphin et al., 2016 \cite{Dauphin}). Training such networks is a problem of minimizing a high-dimensional non-convex and  non-smooth objective function, and is often solved by simple first-order methods such as stochastic gradient descent.
Nevertheless, the success of neural network training remains to be understood from a theoretical perspective. Progress has been made in simplified model problems. %Blum \& Rivest (1993) showed that even training a 3-node neural network is NP-hard \cite{Blum}, and 
Shamir (2016) showed learning a simple one-layer fully connected neural network is hard for some specific input distributions \cite{Shamir}. 
Recently, several works (Tian, 2017 \cite{Tian}; Brutzkus \& Globerson, 2017 \cite{Brutzkus}) focused on the geometric properties of loss functions, which is made possible by assuming that the input data distribution is Gaussian. They showed that stochastic gradient descent (SGD) with random or zero  initialization is able to train a no-overlap neural network in polynomial time.\\
Another notable issue is that DNNs contain millions of parameters and lots of redundancies, potentially causing over-fitting and poor generalization \cite{Zhang2016} besides spending  unnecessary computational resources. 
One way to reduce complexity is to sparsify the network weights using an empirical technique called pruning \cite{LeCun1989} so that the non-essential ones are zeroed out with minimal loss of performance \cite{Han2015,Ullrich2017,Molch2017}. Recently a surrogate $\ell_0$ regularization approach based on a continuous relaxation of Bernoulli random variables in the distribution sense is introduced with encouraging results on small size image data sets \cite{Welling}. This motivated our work here to study deterministic regularization of $\ell_0$ via its Moreau envelope and related $\ell_1$ penalties in a one hidden layer convolutional neural network model \cite{Brutzkus}.    
\begin{figure}[ht]
    \centering
    \includegraphics[scale = 0.3]{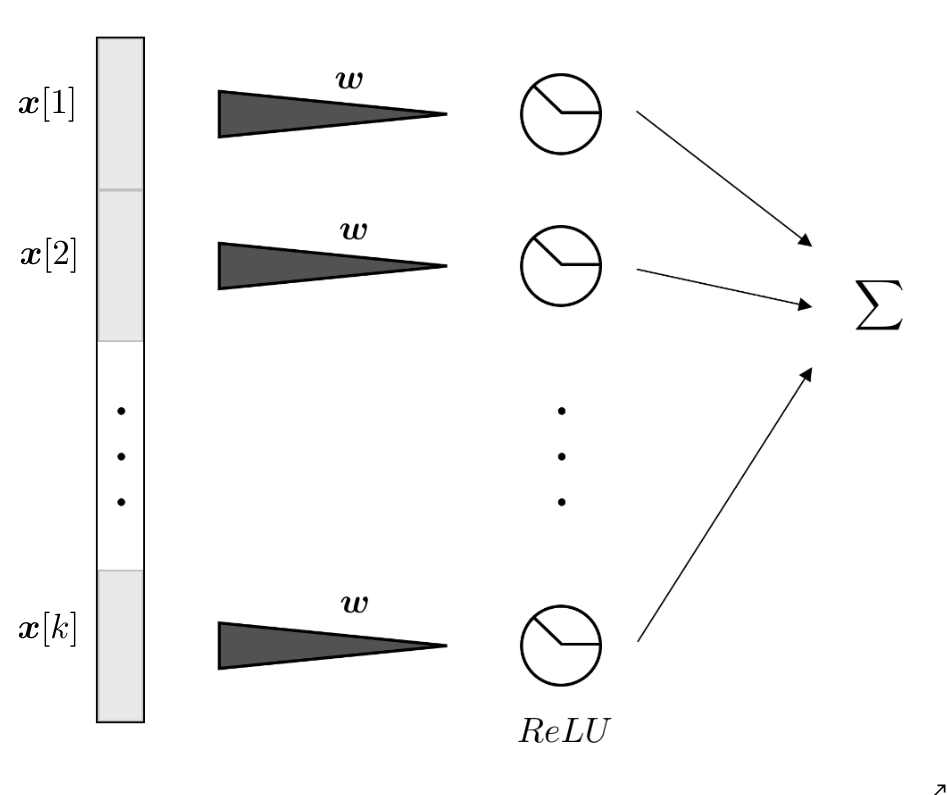}
    \caption{The architecture of a no-overlap neural network}
    \label{network_structure}
\end{figure}

Our contribution: We propose a new method to sparsify DNNs called the Relaxed Variable Splitting Method (RVSM), and prove its convergence on a simple one-layer network (Figure \ref{network_structure}). Consider the population loss:
\begin{equation}
f(\vw) := \mE_{\vx \sim \mathcal{D}} \left[ (L(\vx;\vw)-L(\vx;\vw^*))^2\right].
\end{equation}
where  $L(\vx,\vw)$ is the output of the network with input $\vx$ and weight $\vw$ in the hidden layer. We assume there exists a ground truth $\vw^*$. Consider sparsifying the network by minimizing the Lagrangian
\begin{equation}
\mL_\beta(\vw) = f(\vw) + \|\vw\|_1
\end{equation}
where the $\ell_1$ penalty can be changed to $\ell_0$ or Transformed-$\ell_1$ penalty \cite{Nicolova00,Zhang1}. Empirical experiments show that our method also works on deeper networks, since the sparsification on each layer happens independently of each other.

The rest of the paper is organized as follows. In Section 2, we briefly overview related mathematical results in the study of neural networks and complexity reduction. Preliminaries are in section 3. In Section 4, we state and discuss the main results. The proofs of the main results are in Section 5, and numerical experiments are in Section 6.

\section{Related Work}
In recent years, significant progress has been made in the study of convergence in neural network training. From a theoretical point of view, optimizing (training) neural network is a non-convex non-smooth optimization problem, which is mainly solved by (stochastic) gradient descent. %Blum \& Rivest; Livni et al.; Shalev-Shwartz et al. showed that training a neural network is hard in the worst cases  \cite{Blum,Livni,Shalev}. 
%Shamir showed that if either the target function or input distribution is ``nice", optimization, algorithms used in practice can succeed  \cite{Shamir}.
% Optimization methods in deep neural networks are often categorized into (stochastic) gradient descent methods and others. \\
Stochastic gradient descent methods were first proposed by Robins and Monro in 1951  \cite{Robins}. Rumelhart et al. introduced the popular back-propagation algorithm in 1986  \cite{Rumelhart}. Since then, many well-known SGD methods with adaptive learning rates were proposed and applied in practice, such as the Polyak momentum  \cite{Polyak}, AdaGrad  \cite{ADAGrad}, RMSProp  \cite{RMSProp}, Adam  \cite{Adam}, and AMSGrad  \cite{AMSGrad}.\\ 
The behavior of gradient descent methods in neural networks is better understood when the input has {\it Gaussian} distribution. In 2017, Tian showed the population gradient descent can recover the true weight vector with random initialization for one-layer one-neuron model \cite{Tian}. Brutzkus \& Globerson (2017) showed that a convolution filter with non-overlapping input can be learned in polynomial time \cite{Brutzkus}. Du et al. showed (stochastic) gradient descent with random initialization can learn the convolutional filter in polynomial time and the convergence rate depends on the smoothness of the input distribution and the closeness of patches  \cite{Du1}. Du et al. also analyzed the  polynomial convergence guarantee of randomly initialized gradient descent algorithm for learning a one-hidden-layer convolutional neural network  \cite{Du}.
%A hybrid projected SGD (so called BinaryConnect) is widely used for training various weight quantized DNNs \cite{BC15,YinObj}. 
%Recently, a Moreau envelope based relaxation method (BinaryRelax) is proposed and analyzed to advance 
%weight quantization in DNN training \cite{Yin2018}. 
%Also a blended coarse gradient descent method \cite{YinBlend} is introduced to train fully quantized DNNs in  weights and activation functions, and overcome vanishing gradients.\\
Non-SGD methods for deep learning were also studied in the recent years. Taylor et al. proposed the Alternating Direction Method of Multipliers (ADMM) to transform a fully-connected neural network into an equality-constrained problem to solve  \cite{Taylor}. %Zhang et  al. \cite{Zhang} handled deep supervised hashing (VDSH) problem by an ADMM algorithm to overcome vanishing gradients. Carreira and Wang proposed a method of auxiliary coordinates  (MAC) to replace a nested neural network with a constrained problem without nesting  \cite{Carreira}.\\
A similar algorithm to the one introduced in this paper was discussed in \cite{Lu}. There are a few notable differences. First, their parameter $\varrho$ (respectively our parameter $\beta$) is large (resp. small). Secondly, the update on $\vw$ in our paper does not have the form of an $\text{argmin}$ update, but rather a gradient descent step. Lastly, their analysis does not apply to ReLU neural networks, and the checking step will be costly and impractical for large networks. In this paper, we will show that having $\beta$ small is essential in showing descent of the Lagrangian, angle, and giving a strong error bound on the limit point. We became aware of \cite{Lu} lately after our work was mostly done.
%Wang et al. \cite{Wang}
%proposed a novel Deep Learning Alternating Minimization (DLAM) algorithm to increase computational efficiency and 
%proves global convergence of ADMM for non-convex problems under general assumptions including Lipschitz gradient of the objective functions. 

\section{Preliminaries}
\subsection{The One-layer Non-overlap Network}
In this paper, the input feature $\vx \in \mR^n$ is i.i.d. Gaussian random vector with zero mean and unit variance. Let $\mathcal{G}$ denote this distribution. We assume that there exists a ground truth $\vw^*$ by which the training data is generated. The population risk is then:
\begin{equation}
f(\vw) = \mE_\mathcal{G}[(L(\vx;\vw)-L(\vx;\vw^*))^2].
\end{equation}
We define
\begin{equation}
g(\vu,\vv) = \mE_{\mathcal{G}}[\sigma(\vu \cdot \vx)\sigma(\vv \cdot \vx)].
\end{equation}
Then:
\begin{lem}\cite{Brutzkus,Cho} Assume $\vx \in \mR^d$ is a vector where the entries are i.i.d. Gaussian random variables with mean 0 and variance 1. Given $\vu, \vv \in \mR^d$, denote by $\theta_{\vu,\vv}$ the angle between $\vu$ and $\vv$. Then
\[ g(\vu,\vv) = \frac{1}{2\pi}\|\vu\|\|\vv\|\left(\sin\theta_{\vu,\vv} + (\pi - \theta_{\vu,\vv})\cos\theta_{\vu,\vv}\right). \]	
\end{lem}
For the no-overlap network, the loss function is simplified to:
\begin{equation}\label{simploss}
f(\vw) = \frac{1}{k^2}\big[a(\|\vw\|^2 +  \|\vw^*\|^2) - 2kg(\vw,\vw^*) - 2b\|\vw\|\|\vw^*\|\big].
\end{equation}
where $b = \frac{k^2-k}{2\pi}$ and $a = b + \frac{k}{2}$. 

\subsection{The Relaxed Variables Splitting Method}
Let $\eta > 0$ denote the step size. Consider a simple gradient descent update:
\begin{equation}
\vw^{t+1} = \vw^t - \eta \nabla f(\vw^t).
\end{equation}
It was shown \cite{Brutzkus} that the one-layer non-overlap network can be learned with high probability and in polynomial time. We seek to improve sparsity in the limit weight while also maintain good accuracy. A classical method to accomplish this task is to introduce $\ell_1$ regularization to the population loss function, and the modified gradient update rule. Consider the minimization problem:
\begin{equation}\label{lossFunctionL1}
l(\vw) = f(\vw) + \lambda\|\vw\|_1.
\end{equation}
for some $\lambda > 0$. We propose a new approach to solve this minimization problem, called the Relaxed Variable Splitting Method (RVSM). We first convert (\ref{lossFunctionL1}) into an equation of two variables
\[ l(\vw,\vu) = f(\vw) + \lambda\|\vu\|_1. \]
and consider the augmented Lagrangian
\begin{equation}\label{lagrangian}
\mL_\beta(\vw,\vu) = f(\vw) + \lambda\|\vu\|_1 + \frac{\beta}{2}\|\vw-\vu\|^2.
\end{equation}
Let $S_{\lambda/\beta}(\vw) := sgn(\vw)(|\vw|-\lambda/\beta)\chi_{\{|\vw|>\lambda/\beta\}}$ be the soft thresholding operator. The RSVM is:

% \begin{algorithm}
%     \caption{RVSM Algorithm}\label{RVSMAlgorithm}\begin{algorithmic}[1]
% 	\State\Input The step size $\eta$, parameters $\lambda, \beta$
% 	\State\Initialize $\vw^1,\vu^1$\;
% 	\For{$t=1,2,\ldots,T$}
% 	\State{$\vu^{t+1} \leftarrow \arg\min_{\vu} \mL_\beta(\vw^t,\vu)
% 		=  S_{\lambda/\beta}(\vw^t)$}
% 	\State{$\vw^{t+1} \leftarrow \vw^t - \eta\nabla f(\vw^t) - \eta\beta(\vw^t-\vu^{t+1})$}
%     \EndFor
% 	\State\Output{$\vw^t,\vu^t$}
% \end{algorithmic}
% \end{algorithm}
\begin{algorithm}[H]
		\caption{RVSM
		%Relaxed Variables Splitting Method
		}
		\label{algo:RVSM}\label{RVSMAlgorithm}
		\begin{algorithmic}
			\STATE {\bf Input:} %Hyper-parameters 
			$\eta, \beta,\lambda$, $max_{epoch}$, $max_{batch}$
			\STATE {\bf Initialization:} $\vw^0$
			\STATE {\bf Define:} $\vu^0= S_{\lambda/\beta}(\vw^0)$
			\FOR {$t=0, 1, 2, ..., max_{epoch}$}
			\FOR {$batch=1, 2, ..., max_{batch}$}
			\STATE	$\vw^{t+1} \leftarrow \vw^t - \eta\nabla f(\vw^t) - \eta\beta(\vw^t-\vu^t)$
			\STATE $\vu^{t+1} \leftarrow \arg\min_{\vu} \mathcal{L}_\beta(\vw^t,\vu) = S_{\lambda/\beta}(\vw^t)$
			\ENDFOR
			\ENDFOR 
			
		\end{algorithmic}
	\end{algorithm}

\subsection{Comparison with ADMM}
A well-known, modern method to solve the minimization problem \eqref{lossFunctionL1} is the Alternating Direction Method of Multipliers (or ADMM). In ADMM, we consider the Lagrangian
\begin{equation}\label{lagrangianADMM}
\mL_\beta(\vw,\vu,\vz) = f(\vw) + \lambda\|\vu\|_1 + \langle \vz,\vw-\vu \rangle + \frac{\beta}{2}\|\vw-\vu\|^2.
\end{equation}
and apply the updates:
\begin{equation}\label{eqn ADMM}
\begin{cases}
\vw^{t+1} \leftarrow \arg\min_{\vw} \mL_\beta(\vw,\vu^t,\vz^t)\\
\vu^{t+1} \leftarrow \arg\min_{\vu} \mL_\beta(\vw^{t+1},\vu,\vz^t)\\
\vz^{t+1} \leftarrow \vz^t + \beta(\vw^{t+1}-\vu^{t+1})
\end{cases}
\end{equation}
Although widely used in practice, the ADMM method has several drawbacks when it comes to regularizing deep neural networks: First, the loss function $f$ is often non-convex and only differentiable in some very specific regions, thus the current theory of optimizations does not apply \cite{Wang}. Secondly, 
the update 

\[\vw^{t+1} \leftarrow \arg\min_{\vw} \mL_\beta(\vw^{t+1},\vu,\vz^t)\] 
is not applicable in practice on DNN, as it requires one to know fully how $f(\vw)$ behaves. In most ADMM adaptations on DNN, this step is replaced by a simple gradient descent. Lastly, the Lagrange multiplier $\vz^t$ tends to reduce the sparsity of the limit of $\vu^t$, as it seeks to close the gap between $\vw^t$ and $\vu^t$.
In contrast, the RVSM method  resolves all these difficulties presented by ADMM. First, we will show that in a one-layer non-overlap network, the iterations will keep $\vw^t$ and $\vu^t$ in a nice region, where one can guarantee Lipschitz gradient property for $f(\vw)$. Secondly, 
the update of $\vw^t$ is not an $\arg\min$ update, but rather a gradient descent iteration itself, so our theory does not deviate from practice. Lastly, without the Lagrange multiplier term $\vz^t$, there will be a gap between $\vw^t$ and $\vu^t$ at the limit. The $\vu^t$ is much more sparse than in the case of ADMM, and numerical results showed that $f(\vw^t)$ and $f(\vu^t)$ behave very similarly on deep networks. An intuitive explanation for this is that when the dimension of $\vw^t$ is high, most of its components that will be pruned off to get $\vu^t$ have very small magnitudes, and are often the redundant weights.

In short, the RVSM method is easier to implement (no need to keep track of the variable $\vz^t$), can greatly increase sparsity in the weight variable $\vu^t$, while also maintaining the same performance as ADMM. Moreover, RVSM has convergence guarantee and limit characterization as stated below. 

% \subsection{Parameters Tuning}
% In section \ref{section numerical}, we implement the RVSM on some standard network (ResNet18 and VGG-16) and data sets (CIFAR10 and CIFAR100).  

\section{Main Results}
Before we state our main results, the following Lemma is needed to establish the existence of a Lipschitz constant $L$:
\begin{lem}\label{lemma lipschitz gradient}(Lipschitz gradient)\\
	There exists a global constant $L$ such that the iterations of Algorithm \ref{RVSMAlgorithm} satisfy
	\begin{equation}\label{equation lipschitz gradient}
	\|\nabla f(\vw^t) - \nabla f(\vw^{t+1})\| \leq L \|\vw^t-\vw^{t+1}\|, \;\; \forall t.
	\end{equation}
%	for all $t$.
\end{lem}

An important consequence of Lemma \ref{lemma lipschitz gradient} is: for all $t$, the iterations of Algorithm \ref{RVSMAlgorithm} satisfy:
\[ f(\vw^{t+1})-f(\vw^t) \leq \langle \nabla f(\vw^t), \vw^{t+1}-\vw^t \rangle + \frac{L}{2}\|\vw^{t+1}-\vw^t\|^2. \]

\begin{theorem}\label{RVSMTheorem}
Suppose the initialization of the RVSM Algorithm satisfies:\\
(i) Step size $\eta$ is small so that $\eta \leq \frac{1}{\beta+L}$;\\
(ii) Initial angle $\theta(\vw^0,\vw^*) \leq \pi-\delta$, for some $\delta > 0$;\\
(iii) Parameters $k,\beta,\lambda$ are such that $k \geq 2, \beta \leq \frac{\delta\sin\delta}{k\pi}$, and $\frac{\lambda}{\beta} < \frac{1}{\sqrt{d}}$.\\
Then the Lagrangian $\mL_\beta(\vw^t,\vu^t)$  decreases monotonically; and $(\vw^t,\vu^t)$ converges sub-sequentially to a limit point $(\bar{\vw},\bar{\vu})$, with $\bar{\vu} = S_{\lambda/\beta}(\bar{\vw})$, such that:\\
(i) $0 \in \partial_{\vu} \mL_{\beta}(\bar{\vw},\bar{\vu})$ and $\nabla_{\vw} \mL_{\beta}(\bar{\vw},\bar{\vu})=0$;\\
(ii) $\nabla_{\vw}\,  \mL_{\beta}(\vw^t,\vu^t)=O(\epsilon)$ in $O(1/\epsilon^2)$ iterations;\\
(iii) The limit point $\bar{\vw}$ is close to the ground truth $\vw^*$ in the sense that $\theta(\bar{\vw},\vw^*) < \delta$ and $\|\bar{\vw}-\vw^*\| = O(\beta)$.
\end{theorem}

The full proof of Theorem \ref{RVSMTheorem} is given in the next section. Here we overview the key steps. First, we show that the iterations of Algorithm \ref{RVSMAlgorithm} will eventually bring $\vw^t$ to within a closed annulus $D$ of width $2\,M$ around the sphere centered at origin with radius $\|\vw^*\|$. In other words, there exists a $T$ such that for all $t \geq T, \|\vw^t\| \in [\|\vw^*\|-M,\|\vw^*\|+M]$, for some $0 < M < \|\vw^*\|$. Then with no loss of generality, we can assume that $\vw^t$ is in this closed strip, for all $t$.

Next, for the region $D$ of the iterations, we will show there exists a global constant $L$ such that the Lipschitz gradient property in Lemma \ref{lemma lipschitz gradient} holds.\\
Finally, the Lipschitz gradient property allows us to show the descent of angle $\theta^t$ and Lagrangian   $\mL_{\beta}(\vw^t,\vu^t)$. The conditions on $\eta,\beta,\lambda$ are used to show $\theta^{t+1} \leq \theta^t$; and an analysis of the limit point gives the bound on $\theta(\bar{\vw},\vw^*)$ and $\|\bar{\vw}-\vw^*\|$. From the descent property of $\mL_{\beta}(\vw^t,\vu^t)$, classical results from optimization \cite{Brutzkus} can be used to show that after $T=O\left(\frac{1}{\epsilon^2}\right)$ iterations, we have $\nabla_{\vw} \mL_{\beta}(\vw^t,\vu^t)=O(\epsilon)$, for some $t\in (0, T]$. This leads to the desired convergence rate and finishes the proof. 

It should be noted that without the condition on $\beta$ being small, one may not guarantee monotonicity of $\theta^t$. However, it still can be shown that $\mL_{\beta}(\vw^t,\vu^t)$ decreases and thus the iteration will converge to some limit point $(\bar{\vw},\bar{\vu})$. In this case, the limit point may not be near the ground truth $\vw^*$; i.e. we may not have $\theta(\bar{\vw},\vw^*) < \delta$. Furthermore, the bound on $\|\bar{\vw}-\vw^*\|$ will also be weaker. 

%The proofs of all these key steps do not require convexity of the regularization term $\lambda\|z\|_1$. Thus we can extend our result to other types of regularization penalties. Here we provide another property of the limit point $(\bar{\vz},\bar{\vw})$ and extend our result to $\ell_0$ and transformed $\ell_1$ penalties.

% It follows from Theorem 1 that $\nabla_{\vw} \mL_{\beta}(S_{\lambda/\beta}(\bar{\vw}),\bar{\vw})=\nabla f (\bar{\vw}) + 
% \beta(\bar{\vw}-S_{\lambda/\beta}(\bar{\vw}))=0$. Using the closed form expression 
% of $\nabla f$ (see eq. \eqref{equation gradient} in section \ref{section proof}), we have:

\begin{corollary}\label{Cor thresholded w*}
Suppose the initialization of the RVSM Algorithm satisfies Theorem \ref{RVSMTheorem}, then the $\bar{\vw}$ equation below holds:
%form the limit point $\bar{\vw}$ is some soft-thresholded version of the ground truth $\vw^*$, after some normalization.
%Let $\theta := \theta(\bar{\vw},\vw^*)$, then there exists a constant $0 < C < \frac{1}{1-2k\lambda\sqrt{d}}$ such that
\begin{equation}\label{equation wbar w* normalization}
\vw^* = \frac{k\pi}{\pi-\theta}\beta(\bar{\vw}-S_{\lambda/\beta}(\bar{\vw})) + C\bar{\vw},
\end{equation}

where $\theta := \theta(\bar{\vw},\vw^*)$, 
constant $C\in (0,\frac{1}{1-2k\lambda\sqrt{d}})$.  
Since component-wise, $\bar{\vw}-S_{\lambda/\beta}(\bar{\vw})$ has the same sign as $\bar{\vw}$, the ground truth $\vw^*$ is an expansion of $C\, \bar{\vw}$, or equivalently $\bar{\vw}$ is (up to scalar multiple) a  shrinkage of $\vw^*$. 
\end{corollary}
% This explains theoretically the sparsity in $\bar{\vw}$. 
The proofs of Theorem 1 and Corollary 1.1 do not require convexity of the regularization term $\lambda\|\vu\|_1$, hence extend to other sparse penalties such as $\ell_0$ and transformed $\ell_1$ penalty \cite{Zhang1}. We have:

% \begin{definition}\label{TL1Definition}
% The transformed $l_1$ (TL1) function $\rho_a(x)$ is defined as
% \[ \rho_a(x) = \frac{(a+1)|x|}{a+|x|}, \]
% for some parameter $a \in (0,+\infty)$. With the change of parameter $'a'$, TL1 interpolates $l_0$ and $l_1$ norms:
% \[ \lim_{a \to 0^+} \rho_a(x) = I_{\{x \ne 0\}}, \quad \lim_{a \to +\infty} \rho_a(x) = |x|. \]
% For a vector $x \in \mR^d$, we define
% \[ P_a(x) = \sum_{i=1}^d \rho_a(x_i). \]
% \end{definition}
% More details and discussions of the transformed $\ell_1$ norm are provided the Appendix.

\begin{corollary}\label{Cor}
Under the conditions of Theorem 1 however with the $l_1$ penalty replaced by an $\ell_0$ or transformed-$\ell_1$ penalty, the RVSM Algorithm  converges sub-sequentially to a limit point %when the $l_1$ penalty term is replaced with an $\ell_0$ or TL1 penalty. 
%The limit point 
$(\bar{\vw},\bar{\vu})$ satisfying $\nabla_{\vw} \mL_{\beta}(\bar{\vw},\bar{\vu})=0$. The Lagrangian and angle $\theta^t$ also decrease monotonically, with the limit angle satisfying $\theta(\bar{\vw},\vw^*) < \delta$. Here $\bar{\vu}$ is a thresholding  of $\bar{\vw}$, and equation (\ref{equation wbar w* normalization}) holds with $S_{\lambda/\beta}(\cdot)$ replaced by the thresholding operator of the corresponding  penalty.
\end{corollary}

\section{Proof of Main Results}\label{section proof}
The following Lemmas will be needed to prove Theorem \ref{RVSMTheorem}:
\begin{lem}\label{lemma properties of gradient}(Properties of the gradient, \cite{Brutzkus})\\
For the loss function $f(\vw)$ of equation (\ref{simploss}), the following holds:\\
1. $f(\vw)$ is differentiable if and only if $\vw \ne 0$.\\
2. For $k>1, f(\vw)$ has three critical points:\\
(a) A local maximum at $\vw=0$.\\
(b) A unique global minimum at $\vw = \vw^*$.\\
(c) A degenerate saddle point at $\vw = -\left(\frac{k^2-k}{k^2+(\pi-1)k}\right)\vw^*$.\\
For $k=1, w=0$ is not a local maximum and the unique global minimum $\vw^*$ is the only differentiable critical point.\\ Given $\theta := \theta(\vw,\vw^*)$, the gradient of $f$ can be expressed as
\begin{align}\label{equation gradient}
\begin{split}
\nabla f(\vw) = &\frac{1}{k^2}
\bigg[ 
\bigg( k + \frac{k^2-k}{\pi} - \frac{k}{\pi}\frac{\|\vw^*\|}{\|\vw\|}\sin\theta 
- \frac{k^2-k}{\pi}\frac{\|\vw^*\|}{\|\vw\|} \bigg) \vw - \frac{k}{\pi}(\pi-\theta)\vw^*
\bigg].
\end{split}
\end{align}
\end{lem}

\begin{lem}\label{lemma lipschitz gradient coplanar}(Lipschitz gradient with co-planar assumption, \cite{Brutzkus})\\
	Assume $\|\vw_1\|,\|\vw_2\|\geq M$, $\vw_1,\vw_2,\vw^*$ are on the same two dimensional half-plane defined by $\vw^*$, then
	\[ \|\nabla f(\vw_1) - \nabla f(\vw_2)\| \leq L \|\vw_1-w_2\| \]
	for $L=1+\frac{3\|\vw^*\|}{M}$.
\end{lem}

\begin{lem}\label{lemma closed strip norm}
For $k \geq 1$, there exists constants $M_k, T>0$ such that for all $t \geq T$, the iterations of Algorithm \ref{RVSMAlgorithm} satisfy:
\begin{equation}
\|\vw^t\| \in [\|\vw^*\|-M_k,\|\vw^*\|+M_k].
\end{equation}
where $M_k < \|\vw^*\|$, and $M_k \to 0$ as $k \to \infty$.
\end{lem}

From Lemma \ref{lemma closed strip norm}, WLOG, we will assume that $T=0$.

\begin{lem} \label{lemma Descent in w and u} (Descent of $\mL_\beta$ due to $\vw$ update)\\
For $\eta$ small such that $\eta \leq \frac{1}{\beta+L}$, we have
\[ \mL_\beta(\vu^{t+1},\vw^{t+1}) \leq \mL_\beta(\vw^t,\vu^t). \]
\end{lem}

\subsection{Proof of Lemma \ref{lemma lipschitz gradient}}
\begin{figure}[H]
	\begin{tabular}{l r}
		\begin{minipage}[t]{0.56 \linewidth}
			\hspace{.55in}\includegraphics[scale=0.27]{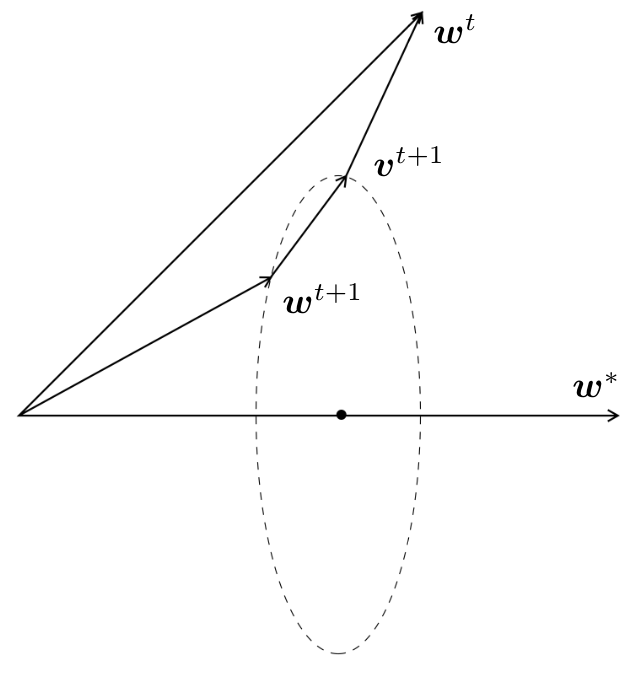} 
			%\caption*{(a) $l_1$ norm} 
		\end{minipage} & 
		\begin{minipage}[t]{0.4\linewidth}
			\includegraphics[scale= 0.27]{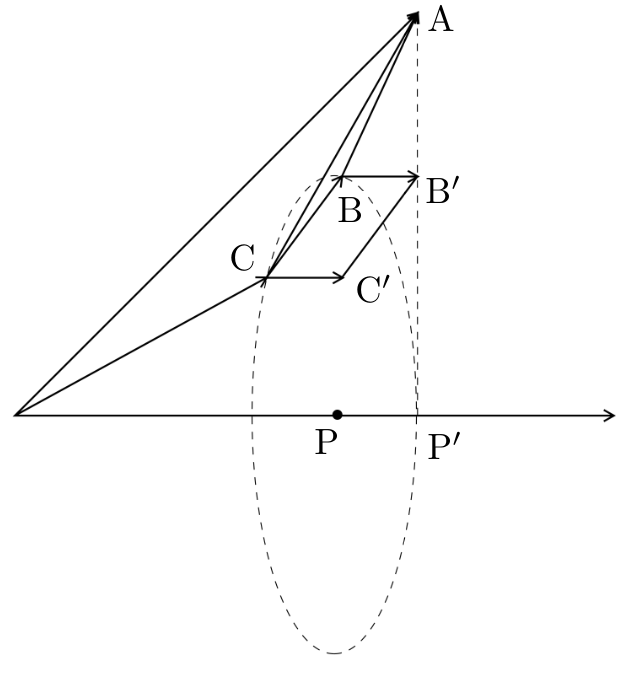} 
			%\caption*{(b) TL1 with $a = 100$} 
		\end{minipage}
	\end{tabular} 
	\caption{Geometry of the update of $\vw^t$ and the corresponding $\vw^{t+1}, \vv^{t+1}$.}
	\label{figure lipschitz triangle}
\end{figure}

% \begin{proof}
By Algorithm \ref{RVSMAlgorithm} and Lemma \ref{lemma closed strip norm}, $\|\vw^t\|\geq \|\vw^*\|-M>0$, for all $t$, and $\vw^{t+1}$ is in some closed neighborhood of $\vw^t$. We consider the subspace spanned by $\vw^t,\vw^{t+1}$, and $\vw^*$, this reduces the problem to a 3-dimensional space.\\
Consider the plane formed by $\vw^t$ and $\vw^*$. Let $\vv^{t+1}$ be the point on this plane, closest to $\vw^t$, such that $\|\vw^{t+1}\| = \|\vv^{t+1}\|$ and $\theta(\vw^{t+1},\vw^*) = \theta(\vv^{t+1},\vw^*)$. In other words, $\vv^{t+1}$ is the intersection of the plane formed by $\vw^t,\vw^*$ and the cone with tip at zero, side length $\|\vw^{t+1}\|$, and main axis $\vw^*$ (See Figure \ref{figure lipschitz triangle}). Then

\begin{align*}
&\|\nabla f(\vw^t) - \nabla f(\vw^{t+1})\|\\ 
\leq &
\|\nabla f(\vw^t) - \nabla f(\vv^{t+1})\| + \|\nabla f(\vv^{t+1}) - \nabla f(\vw^{t+1})\|\\
\leq & L_1\|\vw^t-\vv^{t+1}\| + L_2\|\vv^{t+1}-\vw^{t+1}\| \numberthis \label{equation lipschitz gradient triangle}
\end{align*}
for some constants $L_1, L_2$. The first term is bounded since $\vw^t,\vv^{t+1},\vw^*$ are co-planar by construction, and Lemma \ref{lemma lipschitz gradient coplanar} applies. The second term is bounded by applying Equation \ref{equation gradient} with $\|\vw^{t+1}\| = \|\vv^{t+1}\|$ and $\theta(\vw^{t+1},\vw^*) = \theta(\vv^{t+1},\vw^*)$. It remains to show there exists a constant $L_3>0$ such that
\[ \|\vw^t-\vv^{t+1}\| + \|\vv^{t+1}-\vw^{t+1}\| \leq L_3\|\vw^t-\vw^{t+1}\| \]
Let $A,B,C$ be the tips of $\vw^t,\vv^{t+1},\vw^{t+1}$, respectively. Let $P$ be the point on $\vw^*$ that is at the base of the cone (so $P$ is the center of the circle with $B, C$ on the arc). We will show there exists a constant $L_3$ such that
\begin{equation}\label{equation triangle}
|AB|+|BC| \leq L_3|AC|
\end{equation}
\underline{Case 1:} $A,B,P$ are collinear: By looking at the cross-section of the plane formed by $AB,AC$, it can be seen that $AC$ is not the smallest edge in $\triangle ABC$. Thus there exists some $L_3$ such that Equation \ref{equation triangle} holds.\\
\underline{Case 2:} $A,B,P$ are not collinear: Translate $B,C,P$ to $B',C',P'$ such that $A,B',P'$ are collinear and $BB', CC', PP' \parallelsum \vw^*$. Then by Case 1:
\[|AB'|+|B'C'| \leq L_3|AC'|\]
and $AC'$ is not the smallest edge in $\triangle AB'C'$. By back-translating $B',C'$ to $B,C$, it can be seen that $AC$ is again not the smallest edge in $\triangle ABC$. This implies 
\[|AB|+|BC| \leq L_4|AC|\]

for some constant $L_4$. Thus Equation \ref{equation triangle} is proved. Combining with Equation \ref{equation lipschitz gradient triangle}, Lemma \ref{lemma lipschitz gradient} is proved.
% \end{proof}
\subsection{Proof of Lemma \ref{lemma closed strip norm}}
% \begin{proof}
First we show that if $\|\vw^t\|<\|\vw^*\|$, then the update of Algorithm \ref{RVSMAlgorithm} will satisfy $\|\vw^{t+1}\| > \|\vw^t\|$. By Lemma \ref{lemma properties of gradient}, 
\begin{align*}
\nabla f(\vw) 
= &\frac{1}{k^2}
\bigg[ 
\bigg( k + \frac{k^2-k}{\pi} - \frac{k}{\pi}\frac{\|\vw^*\|}{\|\vw\|}\sin\theta 
- \frac{k^2-k}{\pi}\frac{\|\vw^*\|}{\|\vw\|} \bigg) \vw - \frac{k}{\pi}(\pi-\theta)\vw^*
\bigg]\\
=& \frac{1}{k^2}(C_1 \vw - C_2\vw^*)
\end{align*}
so the update of $\vw^t$ reads
\[ \vw^{t+1} = \vw^t -\eta\frac{C_1^t+\beta k^2}{k^2}\vw^t + \eta \frac{C_2^t}{k^2}\vw^* + \eta\beta \vu^{t+1}, \]
where $C_2^t > 0$. Since $\vu^{t+1}=S_{\lambda/\beta}(\vw^t)$, the term $\eta\beta \vu^{t+1}$ will increase the norm of $\vw^t$. For the remaining terms,
\begin{align*}
C_1^t &= k + \frac{k^2-k}{\pi} - \frac{k}{\pi}\frac{\|\vw^*\|}{\|\vw^t\|}\sin\theta
- \frac{k^2-k}{\pi}\frac{\|\vw^*\|}{\|\vw^t\|}\\
&\leq k + \frac{k^2-k}{\pi} \left(1-\frac{\|\vw^*\|}{\|\vw^t\|}\right)
\end{align*}
When $\frac{\|\vw^*\|}{\|\vw^t\|}$ is large, $C_1^t$ is negative. The update will increase the norm of $\|\vw^t\|$ if $C_1^t+\beta k^2 \leq 0$ and
\[\left\|\frac{C_1^t+\beta k^2}{k^2}\vw^t\right\| > \left\| \frac{C_2^t}{k^2}\vw^* \right\|\]
This condition is satisfied when
\[-\left[ k + \frac{k^2-k}{\pi} \left(1-\frac{\|\vw^*\|}{\|\vw^t\|}\right) + \beta k^2 \right] > 
 \frac{k}{\pi}\frac{\|\vw^*\|}{\|\vw^t\|} \]
When $\frac{\|\vw^*\|}{\|\vw^t\|}>1$, the LHS is $O(k^2)$, while the RHS is $O(k)$. Thus there exists some $M_k$ such that $\vw^t$ will eventually stay in the region $\|\vw^t\| \geq \|\vw^*\|-M_k$. Moreover, as $k \to \infty$, we have $M_k \to 0$.\\
Next, when $\|\vw^t\|>\|\vw^*\|$, the update of $\vw^t$ reads
\[ \vw^{t+1} = \vw^t - \eta\frac{C_1^t}{k^2}\vw^t + \eta\frac{C_2^t}{k^2}\vw^* - \eta\beta(\vw^t-\vu^{t+1}) \]
the last term decreases the norm of $\vw^t$. In this case, $C_1^t$ is positive and
\[ C_1^t \geq \frac{k\pi-k}{\pi} + \frac{k^2-k}{\pi} \left(1-\frac{\|\vw^*\|}{\|\vw^t\|}\right) \]
The update will decrease the norm of $\vw^t$ if
\[ \frac{k\pi-k}{\pi} + \frac{k^2-k}{\pi} \left(1-\frac{\|\vw^*\|}{\|\vw^t\|}\right)
>  \frac{k}{\pi}\frac{\|\vw^*\|}{\|\vw^t\|} \]
which holds when $\frac{\|\vw^*\|}{\|\vw^t\|} < 1$, and the Lemma is proved.
% \end{proof}

\subsection{Proof of Lemma \ref{lemma Descent in w and u}}
\begin{proof}
By the update of $\vu^t$, $\mL_\beta(\vw^t,\vu^{t+1}) \leq \mL_\beta(\vw^t,\vu^t)$. For the update of $\vw^t$, notice that
% \[ \vw^{t+1} = \vw^t - \eta\nabla f(\vw^t) - \eta\beta(\vw^t-\vu^{t+1})\] 
\[ \nabla f(\vw^t) = \frac{1}{\eta}\left(\vw^t-\vw^{t+1}\right) - \beta(\vw^t-\vu^{t+1}) \]
Then for a fixed $\vu:=\vu^{t+1}$, we have
\begin{align*}
&\mL_\beta(\vw^{t+1},\vu) - \mL_\beta(\vw^t,\vu)\\
=& f(\vw^{t+1}) - f(\vw^t) 
+ \frac{\beta}{2}\left(\|\vw^{t+1}-\vu\|^2-\|\vw^t-\vu\|^2\right)\\
\leq& \langle \nabla f(\vw^t), \vw^{t+1}-\vw^t \rangle + \frac{L}{2}\|\vw^{t+1}-\vw^t\|^2\\
+& \frac{\beta}{2}\left(\|\vw^{t+1}-\vu\|^2-\|\vw^t-\vu\|^2\right)\\
=&\frac{1}{\eta}\langle \vw^t - \vw^{t+1}, \vw^{t+1}-\vw^t \rangle
- \beta \langle \vw^t-\vu	, \vw^{t+1}-\vw^t \rangle \\
+& \frac{L}{2}\|\vw^{t+1}-\vw^t\|^2
+ \frac{\beta}{2}\left(\|\vw^{t+1}-\vu\|^2-\|\vw^t-\vu\|^2\right)\\
=& \frac{1}{\eta}\langle \vw^t - \vw^{t+1}, \vw^{t+1}-\vw^t \rangle
+ \left( \frac{L}{2} + \frac{\beta}{2} \right)\|\vw^{t+1}-\vw^t\|^2 \\ 
+& \frac{\beta}{2}\|\vw^{t+1}-\vu\|^2 - \frac{\beta}{2}\|\vw^t-\vu\|^2 \\
-& \beta \langle \vw^t-\vu	, \vw^{t+1}-\vw^t \rangle
- \frac{\beta}{2} \|\vw^{t+1}-\vw^t\|^2\\
=&\left(\frac{L}{2}+\frac{\beta}{2}-\frac{1}{\eta}\right)\|\vw^{t+1}-\vw^t\|^2
\end{align*}

Therefore, if $\eta$ is small so that $\eta \leq \frac{2}{\beta + L}$, the update on $\vw$ will decrease $\mL_\beta$. 
% \end{proof}

%\begin{lem} \label{Limit of L is stationary}(The limit of $\mL_\beta$ is a stationary point)\\
%When Algorithm \ref{RVSMAlgorithm} converges to $(\vz^*,\vw^*,\vz^*)$, we have $0 \in \nabla \mL_\beta (\vz^*,\vw^*,\vz^*)$.
%\end{lem}
%\begin{proof}
%We have
%\begin{equation}
%\nabla \mL_\beta (\vz^t,\vw^t,\vz^t) = 
%\begin{bmatrix}
%\partial g(\vz^t) - \vz^t \\
%\nabla f(\vw^t) + \vz^t + \beta(\vw^t-\vz^t)\\
%\vw^t-\vz^t
%\end{bmatrix}
%\end{equation}
%where $g(\vz) = \lambda\|z\|_1$ and $\partial g(\vz^t)$ is some sub-gradient of $g$ at $\vz^t$. 
%Recall the update $\vz^{t+1} = \vz^t + \beta(\vz^{t+1}-\vw^{t+1})$. Since $\vz^t \to \vz^*$, this implies $\vw^t-\vz^t \to 0$ and $\vw^* = \vz^*$.\\
%Secondly, the update $\vz^{t+1} = \arg\min_u \mL_\beta(z,\vw^t,\vz^t)$ implies $0 \in \partial g(\vz^t)-\vz^t$.\\
%Finally, for $\vw^t$, we have the update $\vw^{t+1} = \vw^t-\eta \nabla f(\vw^t) - \eta z - \eta\beta(\vw^t-\vz^{t+1})$. Since $\vw^t \to \vw^*$ and $\vz^t \to \vz^*$, it follows that 
%\[ \nabla f(\vw^t) + \vz^t + \beta(\vw^t-\vz^t) \to \nabla f(\vw^t) + \vz^t + \beta(\vw^t-\vz^{t+1}) \to 0 \]
%As each component of $\nabla \mL_\beta (\vz^t,\vw^t,\vz^t)$ converges to 0, we have $0 \in \nabla \mL_\beta (\vz^*,\vw^*,\vz^*)$.
%\end{proof}

\subsection{Proof of Theorem \ref{RVSMTheorem}}
% \begin{proof}
We will first show that if $\theta(\vw^0,\vw^*) \leq \pi-\delta$, then $\theta(\vw^t,\vw^*) \leq \pi-\delta$, for all $t$. We will show $\theta(\vw^1,\vw^*) \leq \pi-\delta$, the statement is then followed by induction. To this end, by the update of $\vw^t$, one has
\begin{align*}
% \vw^1 &= C\vw^0 + \left(\eta\frac{\pi-\theta(\vw^0,\vw^*)}{k\pi}\right)\vw^* + \eta\beta \vu^1\\
&= C\vw^0 + \eta \frac{\pi-\theta(\vw^0,\vw^*)}{k\pi}\vw^* + \eta\beta \vu^1
\end{align*}
for some constant $C > 0$. Since $\vu^1 = S_{\lambda/\beta}(\vw^0), \theta(\vu^1,\vw^0) \leq \frac{\pi}{2}$. Notice that the sum of the first two terms on the RHS brings the vector closer to $\vw^*$, while the last term may behave unexpectedly. Consider the worst case scenario: $\vw^0,\vw^*,\vu^1$ are co-planar with $\theta(\vu^1,\vw^0) = \frac{\pi}{2}$, and $\vw^*,\vu^1$ are on two sides of $\vw^0$ (See Figure \ref{worst case}). 
\begin{figure}
\centering
\includegraphics[scale = 0.5]{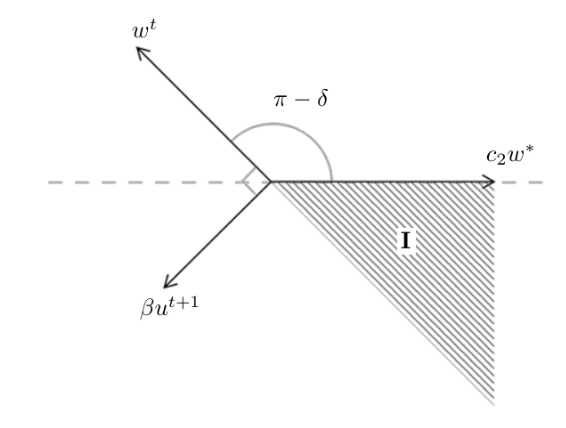}
\caption{Worst case of the update on $\vw^t$}
\label{worst case}
\end{figure}
We need $\frac{\delta}{k\pi} \vw^* + \beta \vu^1$ to be in region I. This condition is satisfied when $\beta$ is small such that
\[ \sin\delta \geq \frac{\beta\|\vu^1\|}{\frac{\delta}{k\pi}\|\vw^*\|}
= \frac{k\pi\beta\|\vu^1\|}{\delta} \]
since $\|\vu^1\| \leq 1$, it is sufficient to have
$\beta \leq \frac{\delta\sin\delta}{k\pi}$.

Next, consider the limit of the RVSM algorithm. Since $\mathcal{L}_\beta(\vw^t,\vu^t)$ is non-negative, by Lemma \ref{lemma Descent in w and u}, $\mathcal{L}_\beta$ converges to some limit $\mathcal{L}$. This implies $(\vw^t,\vu^t)$ converges to some stationary point $(\bar{\vw},\bar{\vu})$. By Lemma \ref{lemma properties of gradient} and the update of $\vw^t$, we have

\begin{equation}\label{EquilibriumCondition}
\overline{\vw} = c_1 \overline{\vw} + \eta c_2\vw^* + \eta\beta \overline{\vu}
\end{equation}
for some constant $c_1> 0, c_2 \geq 0$,
where $c_2 = \frac{\pi - \theta}{k\pi}$, with $\theta := \theta(\bar{\vw},\vw^*)$, and $\bar{\vu} = S_{\lambda/\beta} (\bar{\vw})$. If $c_2=0$, then we must have $\bar{\vw} \parallelsum \bar{\vu}$. But since $\bar{\vu} = S_{\lambda/\beta}$, this implies all non-zero components of $\bar{\vw}$ are either equal in magnitude, or all have magnitude smaller than $\frac{\lambda}{\beta}$. The latter case is not possible when $\frac{\lambda}{\beta} < \frac{1}{\sqrt{d}}$. Furthermore, $c_2=0$ when $\theta(\bar{\vw},\vw^*)=\pi$ or $0$. We have shown that $\theta(\bar{\vw},\vw^*) \leq \pi-\delta$, thus $\theta(\bar{\vw},\vw^*)=0$. Thus, $\bar{\vw} = \vw^*$, and all non-zero components of $\vw^*$ are equal in magnitude. This has probability zero if we assume $\vw^*$ is initiated uniformly on the unit circle. Hence we will assume that almost surely, $c_2>0$. Expression (\ref{EquilibriumCondition}) implies
\begin{equation}\label{ParallelCondition}
c_2 \vw^* + \beta\bar{\vu} \parallelsum \bar{\vw}
\end{equation}
Expression (\ref{ParallelCondition}) implies $\bar{\vw},\bar{\vu}$, and $\vw^*$ are co-planar. Let $\gamma := \theta(\bar{\vw},\bar{\vu})$. From expression (\ref{ParallelCondition}), and the assumption that $\|\vw^*\|=1$, we have
\begin{equation*}
(\langle c_2\vw^*+ \beta\bar{\vu}, \bar{\vw} \rangle)^2 = \|c_2\vw^*+\beta\bar{\vu}\|^2\|\bar{\vw}\|^2
\end{equation*}
or
\begin{align*}
&\|\bar{\vw}\|^2(c_2^2\cos^2\theta + 2c_2\beta\|\bar{\vu}\|\cos\theta\cos\gamma + \beta^2\|\bar{\vu}\|^2\cos^2\gamma)\\
= &\|\bar{\vw}\|^2(c_2^2 + 2c_2\beta\|\bar{\vu}\|\cos(\theta+\gamma) + \beta^2\|\bar{\vu}\|^2)
\end{align*}
This reduces to
\begin{align*}
c_2^2\sin^2\theta - 2c_2\beta\|\bar{\vu}\|\sin\theta\sin\gamma + \beta^2\|\bar{\vu}\|^2\sin^2\gamma = 0,
\end{align*}
which implies
% \begin{equation}\label{Sine Relation}
$\frac{\pi-\theta}{k\pi}\sin\theta = \beta\|\bar{\vu}\|\sin\gamma$.
% \end{equation}
%For small $\beta$ and $\lambda \ll \beta$, the RHS is close to zero, which implies $\theta$ is very close to zero or $\pi$. Since $\theta \leq \pi-\delta$, $\theta$ must be very close to zero.
By the initialization $\beta \leq \frac{\delta\sin\delta}{k\pi}$, we have
$\frac{\pi-\theta}{k\pi}\sin\theta < \frac{\delta}{k\pi}\sin\delta$. This implies $\theta < \delta$.\\
Finally, the limit point satisfies $\|\nabla f(\bar{\vw}) + \beta(\bar{\vw}-\bar{\vu})\| = 0$. By the initialization requirement, we have $\|\beta(\bar{\vw}-\bar{\vu})\| < \beta \leq \frac{\delta\sin\delta}{k\pi}$. This implies $\|\nabla f(\bar{\vw})\| \leq \frac{\delta\sin\delta}{k\pi}$. By the Lipschitz gardient property in Lemma \ref{lemma lipschitz gradient} and critical points property in Lemma \ref{lemma properties of gradient}, $\bar{\vw}$ must be close to $\vw^*$. In other words, $\|\bar{\vw}-\vw^*\|$ is comparable to the chord length of the circle of radius $\|\vw^*\|$ and angle $\theta$:
\begin{align*} \|\bar{\vw}-\vw^*\| &= O\left(2\sin\left(\frac{\theta}{2}\right)\right) = O(\sin\theta)\\
&= O\left(\frac{k\pi\beta\|\bar{\vu}\|\sin\gamma}{\pi-\theta}\right) = O(k\beta\sin\gamma). 
\end{align*}
\end{proof}

\section{Numerical Experiments}\label{section numerical}
First, we experiment RVSM with VGG-16 on the CIFAR10 data set. Table \ref{table vgg} shows the result of RVSM under different penalties. The parameters used are $\lambda = 1.e-5, \beta = 1.e-2$, and $a=1$ for T$\ell_1$ penalty. It can be seen that RVSM can maintain very good accuracy while also promotes good sparsity in the trained network. Between the penalties, $\ell_0$ gives the best sparsity, $\ell_1$ the best accuracy, and T$\ell_1$ gives a middle ground between $\ell_0$ and $\ell_1$. Since the only difference between these parameters is in the pruning threshold, in practice, one may simply stick to $\ell_0$ regularization and just fine-tune the hyper-parameters.

Secondly, we experiment our method on ResNet18 and the CIFAR10 data set. The results are displayed in Table \ref{table resnet18}. The base model was trained on 200 epochs using standard SGD method with initial learning rate 0.1, which decays by a factor of 10 at the 80th, 120th, and 160th epochs. For the RVSM method, we use $\ell_0$ regularization and set $\lambda =$ 1.e-6, $\beta =$ 8.e-2. For ADMM, we set the pruning threshold to be 60\% and $\rho=$1.e-2. The ADMM method implemented here is per \cite{ADMM}, an ``empirical variation'' of the true ADMM (Eq. \ref{eqn ADMM}). In particular, the $\arg\min$ update of $\vw^t$ is replaced by a gradient descent step. Such ``modified'' ADMM is commonly used in practice on DNN.
%\medskip

It can be seen in Table \ref{table resnet18} that RVSM runs quite effectively on the benchmark deep network, promote much better sparsity than ADMM (93.70\% vs. 47.08\%), and has slightly better performance. The  sparsity here is the percentage of zero components over all network  weights.

\begin{table}[H]
  \centering
  \caption{Sparsity and accuracy of RVSM under different penalties on VGG-16 on CIFAR10.}
  \medskip
  
    \begin{tabular}{|l|rr|}
    \hline
    Penalty & Accuracy & Sparsity \\
    \hline
    Base model & 93.82 & 0 \\
    $\ell_1$ & 93.7 & 35.68 \\
    T$\ell_1$  & 93.07 & 63.34  \\
    $\ell_0$ & 92.54 & 86.89 \\
    \hline
    \end{tabular}%
  \label{table vgg}%
\end{table}% 

\begin{table}[H]
  \centering
  \caption{Comparison between ADMM and RVSM ($\ell_0$) for ResNet18 training on the CIFAR10 dataset.}
  \medskip
  
    \begin{tabular}{|l|rr|}
    \hline
    ResNet18     & Accuracy    & Sparsity   \\
    \hline
    SGD & 95.07 & 0 \\
    ADMM  & 94.84 & 47.08  \\
    RVSM ($\ell_0$) & 94.89 & 93.7$0$ \\
    \hline
    \end{tabular}%
  \label{table resnet18}%
\end{table}%

% The channel norms of ResNet18 under these 2 methods are shown in Figure \ref{fig: histogram}. The ADMM method has many channels with very small norms, but not as many with zero norms. This phenomenon was discussed in the Preliminaries, where the Lagrange multiplier term $\vz^t$ in ADMM reduces the sparsity of the weight vectors.

% \begin{figure}[!ht]
% \centering
% \begin{tabular}{cc}
% \includegraphics[width=0.47\columnwidth]{rvsm1.png}&
% % \includegraphics[width=0.45\columnwidth]{rvsm2.png}\\
% \includegraphics[width=0.47\columnwidth]{admm1.png}&
% % \includegraphics[width=0.45\columnwidth]{admm2.png}\\
% (a)  & (b)  \\
% \end{tabular}
% \caption{Channel norms of the adversarially trained ResNet20 under RVSM and ADMM iterations.}
% \label{fig: histogram}
% \end{figure}

\section{Conclusion}
We proved the global convergence of RVSM to sparsify a convolutional ReLU network on a regression problem and analyzed the sparsity of the limiting weight vector as well as its error estimate from the ground truth (i.e. the global minimum). The proof used geometric argument to establish angle and Lagrangian descent properties of the iterations thereby overcame the non-existence  of gradient at the origin of the loss function. Our experimental results provided additional support for the effectiveness of
RVSM via $\ell_0$, $\ell_1$ and T$\ell_1$ penalties on standard deep networks and CIFAR-10 image data. In future work, we plan to extend RVSM theory to multi-layer network and structured (channel/filter/etc.) pruning.

 \section{Acknowledgments}
 The authors would like to thank Dr. Yingyong Qi for suggesting reference \cite{Welling}, and Dr. Penghang Yin for helpful discussions. The work was partially supported by NSF grants DMS-1522383,  IIS-1632935, and DMS-1854434.
% \medskip
%\clearpage

\bibliography{tdbib}
\bibliographystyle{spmpsci}

\end{document}